\documentclass[11pt, leqno, a4paper]{amsart}
\usepackage{amsmath,amstext,amsthm,amsfonts}
\usepackage{amssymb,latexsym}
\usepackage[cp1252]{inputenc}
\usepackage[T1]{fontenc}
\usepackage{lmodern}
\usepackage{hyperref}
\usepackage[english]{babel}

\newtheorem{theorem}{\textbf{Theorem}}[section]
\newtheorem{proposition}[theorem]{\textbf{Proposition}}
\newtheorem{prop}[theorem]{\textbf{Proposition}}

\newtheorem{lemma}[theorem]{\textbf{Lemma}}
\newtheorem{lem}[theorem]{\textbf{Lemma}}
\newtheorem{corollary}[theorem]{Corollary}

\theoremstyle{definition}
\theoremstyle{remark}

\def\pt{\partial_t}

\newcommand{\field}[1]{\mathbb{#1}}
\newcommand{\real}{\field{R}} 

\newcommand{\be} {\beta}        
    
       \newcommand{\De}{\Delta}

\newcommand{\ka} {\kappa}
\newcommand{\la} {\lambda}

\newcommand{\si} {\sigma}

\newcommand{\resp}{\emph{resp. }}

\newcommand{\HH}{\mathbb{H}}
\newcommand{\R}{\mathbb{R}}
\newcommand{\PP}{\mathbb{P}}

\newcommand{\Mh}{\widehat{M}}
\newcommand{\gh}{\hat{g}}
\newcommand{\bh}{\hat{b}}

\newcommand{\Fh}{\widehat{F}}
\newcommand{\pf}{\par{\noindent\textbf{Proof.~}}}

\begin{document}


\title[Eigenvalue estimates for hypersurfaces in $\HH^m \times \R$]
{Eigenvalue estimates for hypersurfaces in $\HH^m \times \R$ and
applications}

\author[B\'erard, Castillon, Cavalcante]{Pierre B\'erard,  Philippe Castillon,
Marcos Cavalcante}

\date{July 6, 2010}
\maketitle

\thispagestyle{empty}


\begin{abstract}
\noindent In this paper, we give a lower bound for the spectrum of
the Laplacian on minimal hypersurfaces immersed into $\HH^m \times
\R$. As an application, in dimension 2, we prove that a complete
minimal surface with finite total extrinsic curvature has finite
index. On the other hand, for stable, minimal surfaces in $\HH^3$ or
in $\HH^2 \times \R$, we give an upper bound on the infimum of the
spectrum of the Laplacian and on the volume growth.
\end{abstract}\bigskip

\bigskip



\bigskip

\section{Introduction}\label{s-in}

In this paper we give a lower bound on the infimum of the spectrum
of the Laplacian $\Delta_g$ on a complete, orientable hypersurface
$(M^m,g)$ minimally immersed into $(\HH^m \times \R, \gh)$ equiped
with the product metric,  with an application to the finiteness of
the index in dimension $2$. In  dimension $2$, under the assumption
that the minimal surface is stable, we give an upper bound on the
infimum of the spectrum and on the volume growth. We also consider
the case when the minimal surface has finite index.

\bigskip

Let us fix some notations.
 Let $\nu$ denote a unit normal field along
$M$ and let $v = \gh(\nu , \partial_t)$ be the component of $\nu$
with respect to the unit vector field $\partial_t$ tangent to the
$\R$-direction in the ambient space.
\bigskip

In Section \ref{s-hr}, we give a lower bound of the spectrum of
$\Delta_g$ which relies on the inequality $- \Delta_g b \ge (m-2) +
v^2$ satisfied by a ``horizontal'' Busemann function $b$ \,(see
Proposition \ref{P-hr1} and Corollary \ref{C-hr2}). In Section
\ref{s-ap}, we give two applications to minimal surfaces in $\HH^2
\times \R$. We prove that a complete minimal surface with finite
total \emph{extrinsic} curvature has finite index (Corollary
\ref{C-hr4}) and we obtain a lower bound for the spectrum of the
Laplacian on a complete  minimal surface contained in a slab
(Proposition \ref{P-slab}). \bigskip

In Section \ref{ss-ir}, we consider the operator $\Delta_g + a +
bK_g$ on a complete Riemannian surface. When $a \geq 0$ and $b> 1/4$,
we show that the positivity of this operator implies an upper bound
on the infimum of the spectrum of $\Delta_g$ and on the volume
growth of $M$ (see Proposition \ref{P-up} and Proposition
\ref{P-vge1}). In Section \ref{ss-aps}, we apply these results to
stable minimal surfaces in $\HH^3$ or $\HH^2 \times \R$,
generalizing and extending results of A. Candel, \cite{C}. Candel
used Pogorelov's method, \cite{P81}. We use the method of Colding
and Minicozzi, \cite{CM02, Cas06}. \bigskip

In Section \ref{s-apf}, we give some applications of our general
lower bounds on the spectrum to higher dimensional hypersurfaces. In
Section \ref{s-pc}, we provide some preliminary technical lemmas.
\bigskip

The third author gratefully acknowledges CAPES and FAPEAL for their
financial support and Institut Fourier for their hospitality during
the preparation of this paper.

\section{Preliminary computations}\label{s-pc}

In this section, we make some preliminary computations for later
reference. For the sake of simplicity, we work in the following
model for the hyperbolic space $\HH^{m+1}$,
\begin{equation}\label{E-pc1}
\left\{%
\begin{array}{ll}
\HH^{m+1}&= \R^m \times \R,\\
h &= e^{2s}\big( dx_1^2 + \cdots + dx_m^2 \big) + ds^2 \text{~ at
the point~} (x,s) \in \HH^{m+1}.
\end{array}%
\right.
\end{equation}

These coordinates are known as ``horocyclic coordinates'' because
the slices $\R^m\times\{s\}$ are horospheres and the coordinate
function $s$ is a Busemann function. They are quite natural when
some Busemann function plays a special role, as will be the case in
the sequel. Let $\gamma_0$ be the geodesic ray
\begin{equation}\label{E-pc2}
\gamma_0~ :~ \left\{%
\begin{array}{l}
[0, \infty) \to \HH^{m+1}, \\
u \mapsto \gamma_0(u) = (0, \ldots, 0,u).
\end{array}%
\right.
\end{equation}

The Busemann function (see \cite{BGS85}, p. 23) associated with
$\gamma_0$ is the function
\begin{equation}\label{E-pc3}
B~ :~ \left\{%
\begin{array}{l}
\HH^{m+1} \to \R, \\
(x,s) \mapsto B(x,s) = s.
\end{array}%
\right.
\end{equation}

In the sequel, we denote by
\begin{equation}\label{E-pc4}
\left\{%
\begin{array}{ll}
D^{h} & \text{~ the Levi-Civita connexion},\\
\Delta_h & \text{~ the geometric (\emph{i.e.} non-negative)
Laplacian},
\end{array}%
\right.
\end{equation}
for the hyperbolic metric $h$ on $\HH^{m+1}$. \bigskip

\begin{lemma}\label{L-pc1}
With the above notations, we have the formulas,
\begin{equation}\label{E-pc5}
\Delta_h B = -m,
\end{equation}
\begin{equation}\label{E-pc6}
\mathrm{Hess}_hB = e^{2s} \big(dx_1^2 + \cdots + dx_m^2 \big)
\end{equation}
at the point $(x,s) \in \HH^{m+1}$. In particular, if we decompose
the vector $u \in T_{(x,s)}\HH^{m+1}$ $h$-orthogonally as $u
=(u_x,u_s)$, we have,
\begin{equation}\label{E-pc7}
\mathrm{Hess}_hB(u,u) = h(u_x,u_x).
\end{equation}
\end{lemma}\bigskip

The proof is straightforward. \qed


\bigskip

Recall the following general lemmas.

\begin{lem}\label{L-pc2}
Let $(M^m,g) \looparrowright (\Mh^{m+1}, \gh)$ be an orientable
isometric immersion with unit normal field $\nu$ and corresponding
normalized mean curvature $H$. Let $\Fh : \Mh \to \R$ be a smooth
function and let $F := \Fh|_{M}$ be its restriction to $M$. Then, on $M$,
$$\Delta_g F = \Delta_{\gh} \Fh|M + \mathrm{Hess}_{\gh}\Fh(\nu,\nu) -
mH d\Fh(\nu).$$
\end{lem}

\pf See for example \cite{CG92}, Lemma 2. \qed


\begin{lemma}\label{L-pc3}
Assume that the manifold $(M,g)$ carries a function $f$ which
satisfies%
$$|df|_g \le 1 \text{~ and ~} - \Delta_g f \geq c \text{~ for some constant ~} c>0.$$
Then, any smooth, relatively compact domain $\Omega \subset M$
satifies the isoperime\-tric inequalities%
$$\mathrm{Vol}_{m-1}(\partial \Omega) \geq c \, \mathrm{Vol}_m(\Omega) \text{~ and ~}
\lambda_1 (\Omega) \geq \frac{c^2}{4},$$ where $\lambda_1(\Omega)$
is the least eigenvalue of $\Delta_g$ in $\Omega$, with Dirichlet
boundary condition.
\end{lemma}

\pf Integration by parts and Cauchy-Schwarz. \qed




\section{Hypersurfaces in $\HH^m \times \R$}\label{s-hr}

We consider orientable, isometric immersions $(M^m,g)
\looparrowright (\Mh^{m+1}, \gh)$, with unit normal $\nu$, where
$\Mh = \HH^m \times \R$ with the product metric $\gh = h + dt^2$. We
take the model (\ref{E-pc1}) for the hyperbolic space (here with
dimension $m$), so that $\Mh$ is the product $\R^{m-1} \times \R
\times \R$, with the Riemannian metric $\gh$ given by $$\gh =
e^{2s}(dx_1^2 + \cdots + dx_{m-1}^2) + ds^2 + dt^2.$$ We define the
function $\bh$ on $\Mh$ by
\begin{equation}\label{E-hr1}
\bh(x_1, \ldots, x_{m-1},s,t) = s.
\end{equation}

This function is in fact a Busemann function of $\Mh$ (seen as a
Cartan-Hadamard manifold) associated with a ``horizontal'' geodesic
(justifying the name ``horizontal'' Busemann function used in the
introduction).\bigskip

We call $b := \bh|_M$ the restriction of $\bh$ to $M$. We decompose
the unit vector $\nu$ according to the product structure $\R^{m-1}
\times \R \times \R$, orthogonally with respect to $\gh$, as
\begin{equation}\label{E-hr2}
\nu = \nu_x + w \partial_s + v \partial_t.
\end{equation}

Applying Lemma \ref{L-pc2}, we obtain the equation
\begin{equation}\label{E-hr3}
\Delta_g b = \Delta_{\gh} \bh|_M + \mathrm{Hess}_{\gh} \bh(\nu,\nu)
- m H \gh(\nu, \partial_s).
\end{equation}

Using (\ref{E-pc7}) and (\ref{E-hr2}), it can we rewritten as
\begin{equation}\label{E-hr4}
- \Delta_g b = (m-1) - |\nu_x|^2 + mH w,
\end{equation}
and we note that $|\nu_x|^2 + v^2 + w^2 = 1$. It follows that
\begin{equation}\label{E-hr5}
- \Delta_g b \geq  (m-2) + v^2 + w^2 - mH |w|.
\end{equation}

For minimal hypersurfaces, we deduce from (\ref{E-hr5}) the
following results.\bigskip

\begin{proposition}\label{P-hr1}
Let $(M^m,g) \looparrowright (\HH^m \times \R, \gh)$ be a complete,
orientable, mi\-nimal hypersurface, with normal vector $\nu$. Recall
that $v = \gh(\nu,\partial_t)$. Then,
    \begin{equation}\label{E-hr6}
    - \Delta_g b \geq (m-2) + v^2.
    \end{equation}
\end{proposition}

\begin{corollary}\label{C-hr2}
Let $(M^m,g) \looparrowright (\HH^m \times \R, \gh)$ be a complete,
orientable, mi\-ni\-mal hypersurface, with normal vector $\nu$. Let
$v = \gh(\nu,\partial_t)$. Let $\lambda_{\sigma}(\Delta_g)$ be the
infimum of the spectrum of the Laplacian $\Delta_g$ on $M$. Then
\begin{equation}\label{E-hr7}
    \lambda_{\sigma}(\Delta_g) \geq \big( \frac{m-2 +
\inf_M v^2}{2} \big)^2 \geq \big( \frac{m-2}{2} \big)^2.
\end{equation}
\end{corollary}

\begin{corollary}\label{C-hr3}
Let $(M^m,g) \looparrowright (\HH^m \times \R, \gh)$ be a complete,
orientable, minimal hypersurface, with $m \ge 3$. Then $(M,g)$ is
non-parabolic.
\end{corollary}

\pf Apply Proposition 10.1 of  \cite{Gri99} using (\ref{E-hr7}). \qed \bigskip

When the mean curvature $H$ is non-zero, we also obtain the
following result from inequality (\ref{E-hr5}),

\begin{proposition}\label{P-hr1a}
Let $(M^m,g) \looparrowright (\HH^m \times \R, \gh)$ be a complete,
orientable hypersurface, with normal vector $\nu$ and constant mean
curvature $H$, $0 \le H \le \frac{m-1}{m}$. Recall that $v =
\gh(\nu,\partial_t)$. Then,%
    \begin{equation}\label{E-hr7a}
    - \Delta_g b \geq (m-2)(1 - \sqrt{1-v^2}) + (m-2)(1 -
    \frac{mH}{m-2}) \sqrt{1-v^2}.
    \end{equation}
\end{proposition}

\textbf{Remarks}. (i) Inequalities (\ref{E-hr6}) and (\ref{E-hr7})
are sharp. Indeed, take the horizontal slice $M = \HH^{m} \times
\{0\}$, in that case $v=1$, or take $M=\PP \times \R$, where $\PP$
is some totally geodesic $(m-1)$-space in $\HH^m$, in that case
$v=0$. (ii) In dimension $2$, Corollary~\ref{C-hr2} is empty in
general. However, inequality (\ref{E-hr6}) is useful even in
dimension $2$, as we will show in Section~\ref{s-ap}. (iii)
Inequality (\ref{E-hr7a}) generalizes an earlier result of the
second author (\cite{Cas97}) for submanifolds immersed in Hadamard
manifolds. We point out that it is more convenient in our context to
use the ``horizontal'' Busemann function rather than the hyperbolic
distance function as in \cite{Cas97}. (iv) The above inequalities
still hold if $M^m$ is only assumed to have mean curvature bounded
from above by $H$.\bigskip


\section{Applications to minimal hypersurfaces in $\HH^m \times
\R$}\label{s-ap}

\subsection{Index of minimal surfaces immersed in $\HH^2 \times
\R$}\label{ss-ap1}

The stability operator of a minimal hypersurface $M^m \looparrowright \HH^m
\times \R$ is given by
\begin{equation}\label{E-ap1}
J_M = \Delta + (m-1)(1-v^2) - |A|^2,
\end{equation}
where $v$ is the vertical component of the unit normal $\nu$, and
$A$ the second fundamental form of the immersion (see \cite{BS08}).
It turns out that the spectrum of the operator $\Delta +
(m-1)(1-v^2)$ is bounded from below by a positive constant. More
precisely, we have the following result.

\begin{proposition}\label{P-hr3}
Let $(M^m,g) \looparrowright (\HH^m \times \R, \gh)$ be a complete,
orientable, minimal hypersurface with normal vector $\nu$. Let $v =
\gh(\nu,\partial_t)$. Then the spectrum of the operator $\Delta_g +
(m-1)(1 - v^2)$ on $M$ is bounded from below by $(\frac{m-1}{2})^2$.
\end{proposition}

\textbf{Proof}. We start from the inequality (\ref{E-hr6}), $-
\Delta_g b \ge (m-2) + v^2$. We multiply this inequality by $f^2$,
where $f \in C^{\infty}_0(M)$, and integrate by parts using the fact
that $|db|_g \le 1$. We obtain (all integrals are taken with respect
to the Riemannian measure $dv_g$),
$$(m-2) \int_M f^2 + \int_M v^2 f^2 \le \int_M |df^2| \le 2 \int_M |f| |df|.$$

We re-write this inequality as
$$(m-1) \int_M f^2 \le 2 \int_M |f| |df| + \int_M (1-v^2) f^2.$$

Using the Cauchy-Schwarz inequality $2|f|.|df| \le \frac{1}{a}|df|^2
+ af^2$ for $a>0$, we obtain
$$a(m-1-a) \int_M f^2 \le \int_M \big( |df|^2 + a(1-v^2) f^2\big) \le
\int_M \big( |df|^2 + (m-1)(1-v^2) f^2\big),$$

provided that $0 \leq a \leq m-1$. We can now maximize the constant in the
left-hand side by choosing $a=(m-1)/2$. \qed\bigskip

\textbf{Remark}. We observe that equality is achieved in the above
inequality when $M$ is a slice $\HH^m \times \{t_0\}$, in which case
$v=1$. If we assume that $v^2 \le \alpha^2 < 1$, the spectrum of
$\Delta_g + (m-1)(1 - v^2)$ is bounded from below by $(m-1)(1 -
\alpha^2)$.

\begin{corollary}\label{C-hr4}
Let $(M^2,g) \looparrowright (\HH^2 \times \R, \gh)$ be a complete,
orientable, minimal surface, with second fundamental form $A$. If
$\int_M |A|^2 dv_g$ is finite, then the immersion has finite index.
\end{corollary}

\textbf{Proof.} When $\int_M |A|^2$ is finite, the second
fundamental form tends to zero uniformly at infinity (see
\cite{BS08}, Theorem 4.1). Using Proposition \ref{P-hr3} with $m=2$,
it follows that the essential spectrum of the Jacobi operator $J_M$
is bounded from below by $\frac{1}{4}$. Since the operator $J_M$ is
also bounded from below, it follows that it has only finitely many
negative eigenvalues (see \cite{BdCS97}, Proposition 1). \qed
\bigskip

\textbf{Remark}. This corollary answers a question raised in
\cite{BS08}, where the finiteness of the index of $J_M$ is proved in
dimension $m \ge 3$ under the assumption that $\int_M |A|^m$ is
finite, and in dimension $2$ under the assumption that both $\int_M
v^2$ and $\int_M |A|^2$ are finite. In dimension $m \ge 3$, the
index of $J_M$ is bounded from above by a constant times $\int_M
|A|^m$ (see \cite{BS08}). In the next section, we investigate bounds
on the index in dimension $2$. \bigskip

\subsection{Bounds on the index of minimal surfaces immersed in $\HH^2 \times
\R$}\label{ss-ap1a}

\begin{prop}\label{P-bd}
Let $(M^2,g) \looparrowright (\HH^2 \times \R, \gh)$ be a complete,
orientable, minimal surface, with second fundamental form $A$. If
$\int_M |A|^2 \, dv_g$ is finite, then for any $r>1$, there exists a
constant $C_r$ such that the index of the immersion is bounded from
above by $C_r \, \int_M |A|^{2r} \, dv_g$.
\end{prop}\bigskip

\textbf{Remarks}. (i) Recall that the assumption that $\int_M |A|^2
\, dv_g$ is finite implies that $A$ tends to zero uniformly at
infinity. It follows that the integrals $\int_M |A|^{2r} \, dv_g$
are all finite. (ii) Our proof provides a constant $C_r$ which tends
to infinity when $r$ tends to $1$. We do not know whether there is a
bound of the index in terms of $\int_M |A|^2 \, dv_g$ as this is the
case for minimal surfaces in $\R^3$ (see \cite{Ty87}). \bigskip

\pf As in Section \ref{ss-ap1}, we write the Jacobi operator as $J =
\Delta_g + 1- v^2 - |A|^2$. The closure $\widetilde{Q}$ of the
quadratic form $Q[f] = \int_M \big( |df|^2  + (1-v^2) f^2 \big) \,
dv_g$ with domain $C_0^1(M)$ satisfies the Beurling-Deny condition
(if $f$ is in the domain of $\widetilde{Q}$, then so is $|f|$ and
$\widetilde{Q}[|f|] = \widetilde{Q}[f]$, see \cite{Da89}, Theorem
1.3.2) and, by Proposition~\ref{P-hr3}, the Cheeger inequality
\begin{equation}\label{A1}
\int_M f^2 \, dv_g \le 4 Q[f], ~~ \forall f \in C_0^1(M).
\end{equation}

On the other-hand, the surface $M$ satisfies the Sobolev inequality
\begin{equation}\label{A2}
\int_M f^2 \, dv_g \le S \big( \int_M |df|_g^2 \, dv_g \big)^2, ~~
\forall f \in C_0^1(M),
\end{equation}
for some constant $S>0$. Indeed, this follows from the Sobolev
inequality for minimal surfaces in $\HH^2 \times \R$, using the fact
that the ambient space has non-positive curvature and infinite
injectivity radius (see \cite{HS74}).
\bigskip

From the above Cheeger and Sobolev inequalities, we can establish
that for any $q \ge 1$, there exists a constant $D_q$ such that for
any $f \in C_0^1(M)$,
\begin{equation}\label{A3}
\big( \int_M |f|^{2q} \, dv_g \big)^{1/q} \le D_q Q[f].
\end{equation}
When $q$ is an integer, the inequality follows from an induction
argument and we can conclude by interpolation. \bigskip

We can then apply Theorem 1.2 of \cite{LS97} to conclude that the
index is less than $e^p D_q^p \int_M |A|^{2p} \, dv_g$ where $p =
q/(q-1)$. \qed \bigskip

\subsection{Hypersurfaces in a slab}\label{ss-ap2}

In this section, we use the computations of Section~\ref{s-hr} to
give a lower bound on the spectrum of the Laplacian on a complete
minimal surface immersed in a slab $\HH^2\times[-a,a]$, $a>0$.
\bigskip

Let us first consider functions on $\HH^m \times \real$ depending
only on the height $t$, namely $\hat{\be}(x,s,t)=f(t)$. In  this
case, $d\hat{\be}=f'(t) dt$, and%
$$
 \mathrm{Hess}_{\gh}\hat{\be}(X,Y)=f''(t)\gh(X,\pt)\gh(Y,\pt).
$$
In particular,%
$$
\De_{\gh}\hat{\be}=-f''(t)\quad\text{ and }\quad
\mathrm{Hess}_{\gh}\hat{\be}(\nu,\nu)=v^2f''(t).
$$

Let us define $\be=\hat{\be}|_M$. Using Lemma \ref{L-pc2},  we have%
\begin{equation}\label{slab1}
-\De_g\be= (1-v^2)f''(t)+mHvf'(t).
\end{equation}
In order to estimate the first eigenvalue of a \emph{minimal}
hypersurface $M^m \looparrowright \HH^m \times \real$, we use the
identity (\ref{slab1}) with some particular choice of $f$. For
instance, let $\hat{\be}(x,s,t)=\frac{1}{2}t^2$. In this case, we
have%
$$
-\De_g \be = (1-v^2).
$$

Assume now that $M^m \looparrowright \HH^m \times [-a,a]$, for some
$a>0$. Then,%
$$
-\De_g \be =  (1-v^2) \text{ ~and~ } |d\be|\leq a.
$$

If we define $Z=b+\be$, where $b$ is the restriction of the Busemann
function $\hat{b}$ to $M^m$, we can use the last inequality in
(\ref{E-hr5}) to obtain
\begin{equation}\label{slab2}
-\De \,Z  \geq   m-1 \text{ ~and~ } |dZ|\leq\sqrt{1+a^2} .
\end{equation}

Using the above notation and Lemma \ref{L-pc3}, we have the
following estimate,

\begin{proposition}\label{P-slab}
Given $a>0$, let $(M^m,g) \looparrowright (\HH^m \times [-a,a],
\gh)$ be a complete, immersed, orientable,  minimal hypersurface.
Then, the infimum of the spectrum of $\Delta_g$ on $M$ is positive.
More precisely,
\begin{equation}\label{slab3}
\la_{\sigma}(\Delta_g) \geq \frac{(m-1)^2}{4(1+a^2)}.
\end{equation}
\end{proposition}

\section{Bounds derived from a stability assumption}\label{s-stab}

Let $(M,g)$ be a complete Riemannian surface with (non-negative)
Laplace operator $\Delta_g$ and Gaussian curvature $K_g$. Let $a, b$
be real numbers, with $a\geq 0$ and  $b > 1/4$. Let $L$ be the operator $L
= \Delta_g + a + b K_g$. \bigskip

Let $\mathrm{Ind}(L,\Omega)$ denote the number of negative
eigenvalues of the operator $L$ in $\Omega$, with Dirichlet boundary
conditions on $\partial \Omega$. The index, $\mathrm{Ind}(L)$, of
the operator $L$ is defined to be the supremum
$$\mathrm{Ind}(L) = \sup \{\mathrm{Ind}(L,\Omega) ~|~ \Omega \Subset M\}$$
taken over the relatively compact subdomains $\Omega$ in
$M$. \bigskip

In Section \ref{ss-ir}, we state two intrinsic consequences of the
assumption that the operator $L$ has finite index. In
Sections~\ref{ss-aps} and \ref{ss-fa}, we consider applications to
minimal and \textsc{cmc} surfaces.
\bigskip

\subsection{Intrinsic results}\label{ss-ir}

\begin{proposition}\label{P-up}
Let $(M,g)$ be a complete non-compact Riemannian surface. Let $a
\geq 0$ and $b > \frac{1}{4}$. Denote by $\Delta_g$ the
(non-negative) Laplacian and by $K_g$ the Gaussian curvature of
$(M,g)$. Denote by $\lambda_{\sigma}(\Delta_g)$ the infimum of the
spectrum of $\Delta_g$ and by $\lambda_{e}(\Delta_g)$ the infimum of
the essential spectrum of $\Delta_g$.
\begin{enumerate}
    \item If the operator $\Delta_g + a + b K_g$ is non-negative on
    $C_0^{\infty}(M)$, then, $$\lambda_{\sigma}(\Delta_g) \le \frac{a}{4b-1}.$$
    \item If the operator $\Delta_g + a + b K_g$ has finite index on
    $C_0^{\infty}(M)$ and if $M$ has infinite volume, then,
    $$\lambda_{e}(\Delta_g) \le \frac{a}{4b-1}.$$
\end{enumerate}
\end{proposition}\bigskip

\pf The proof uses the method of Colding-Minicozzi \cite{CM02}, and
more precisely Lemma 1.8 in the second author's paper \cite{Cas06}.

\emph{Proof of Assertion 1}. We can assume the surface to have
infinite volume (otherwise $\lambda_{\sigma}(\Delta_g) = 0$ because
the function $1$ is in $L^2(M,v_g)$ and the estimate is trivial).
Fix a point $x_0 \in M$ and let $r(x)$ denote the Riemannian
distance to the point $x_0$. Given $S > R > 0$, let $B(R)$ denote
the open geodesic ball in $M$ with center $x_0$ and radius $R$. Let
$C(R,S)$ denote the open annulus $B(S) \setminus \bar{B}(R)$. Let
$V(R)$ denote the volume of $B(R)$ and $L(R)$ the length of its
boundary $\partial B(R)$. Let $G(R)$ denote the integral curvature
of $B(R)$, $G(R) = \int_{B(R)} K_g(x) \, dv_g (x)$, where $dv_g$
denotes the Riemannian measure. The main idea in \cite{Cas06} is to
use the work of Shiohama-Tanaka \cite{ST88, ST93} on the length of
geodesic circles, where it is shown that the function $L(r)$ is
differentiable almost everywhere and related to the Euler
characteristic and to the integral curvature of geodesic balls by
the formula (\cite{Cas06}, Theorem 1.7)
$$L'(r) \le 2 \pi \chi(B(r)) - G(r) \le 2\pi - G(r),$$

where the second inequality comes from the fact that the Euler
characteristic of balls is less than or equal to $1$.
Recall the following lemma.
\bigskip

\begin{lem}[Lemma 1.8. in \cite{Cas06}]\label{L-cas}
For $0 < R < S$, let $\xi : [R,S] \to \R$ be such that $\xi \ge 0,
\xi' \le 0$, $\xi'' \ge 0$ and $\xi(S)=0$. Then
$$
\begin{array}{ll}
\int_{C(R,S)} K_g \xi^2(r) \, dv_g \le & - \xi^2(R) G(R) +  2\pi
\xi^2(R)
- 2 \xi(R) \xi'(R) L(R) \\
& ~~  - \int_{C(R,S)} (\xi^2)''(r) \, dv_g.\\
\end{array}
$$
\end{lem}\bigskip

To prove Assertion 1,we choose $\xi$ as in Lemma \ref{L-cas}, and a
function $f : B(S) \to \R$ such that $f (r) \equiv \xi(R)$ on
$B(R)$, $f (r) = \xi (r)$ on $C(R,S)$, and we write the positivity
assumption,
$$0 \le \int_M |df|_g^2 \, dv_g + a \int_M f^2 \, dv_g
+ b \int_M K_g f^2 \, dv_g.$$

On the ball $B(R)$, we have $$\int_{B(R)} K_g f^2 \, dv_g = \xi^2(R)
G(R) \text{ ~and~ } \int_{B(R)} |df|^2 \, dv_g = 0.$$ Using
Lemma \ref{L-cas}, we obtain
$$
\begin{array}{ll}
0 \le & \int_{C(R,S)} (\xi')^2(r) \, dv_g + a \int_M f^2 \, dv_g
+ b \xi^2(R) G(R) - b \xi^2(R) G(R) \\[4pt]
& + 2 \pi b \xi^2(R)- 2 b \xi(R) \xi'(R) L(R) - b \int_{C(R,S)} (\xi ^2)''(r) \, dv_g, \\
\end{array}
$$
and hence,
\begin{equation}\label{E-star}
\begin{array}{ll}
0 \le & (1-2b) \int_{C(R,S)} (\xi')^2(r) \, dv_g + a \int_M f^2 \, dv_g
\\[4pt]
& + 2 \pi b \xi^2(R) - 2 b \xi(R) \xi'(R) L(R) - 2b \int_{C(R,S)} \xi(r) \xi''(r) \, dv_g. \\
\end{array}
\end{equation}

We choose $\xi(r) = (S-r)^k$ in $[R,S]$ for $k \ge 1$ big enough (we
will eventually let $k$ tend to infinity). Then $\xi(r) \xi''(r) =
(1 - \frac{1}{k}) (\xi'(r))^2$. It follows that
$$
\begin{array}{ll}
0 \le & (1-4b + 2b/k) \int_M |df|^2 \, dv_g + a \int_M f^2 \, dv_g
\\[4pt]
& + 2b \Big( \pi (S-R)^{2k} + kL(R)(S-R)^{2k-1} \Big).
\end{array}
$$

Using the fact that $\int_M f^2 \, dv_g \ge (S-R)^{2k} V(R)$, we
obtain
\begin{equation}\label{E-sstar}
\begin{array}{ll}
\lambda_{\sigma}(\Delta_g) & \le \dfrac{\int_M |df|^2 \,
dv_g}{\int_M f^2 \, dv_g} \\[4pt]
& \le \dfrac{a}{4b-1 - 2b/k} + \dfrac{2b}{(4b-1 - 2b/k)V(R)}\big(
\pi + \frac{k L(R)}{S-R} \big) .\\
\end{array}
\end{equation}

We first let $S$ tend to infinity, then we let $R$ tend to infinity,
using the fact that $M$ has infinite volume, and we let finally $k$
tend to infinity to obtain
$$\lambda_{\sigma}(\Delta_g) \le \dfrac{a}{4b-1}.$$

\emph{Proof of Assertion 2}. It is a well-known fact that the
finiteness of the index of the operator $\Delta_g + a + b K_g$
implies that it is non-negative outside a compact set (see
\cite{FC85}, Proposition 1). We choose $R_0$ big enough for
$\Delta_g + a + b K_g$ to be non-negative in $M \setminus B(R_0)$.
Next, for $S > R > R_1+1
> R_0+1$, we choose $\xi$ as in Lemma \ref{L-cas}, and a test
function $f$ as follows
\begin{equation}\label{E-3star}
f(r) = \left\{%
\begin{array}{ll}
0 & \text{ ~in~} B(R_1),\\
\xi(R)(r-R_1) & \text{ ~in~} C(R_1,R_1+1), \\
\xi(R) & \text{ ~in~} C(R_1+1,R),\\
\xi (r) & \text{ ~in~} C(R,S).\\
\end{array}%
\right.
\end{equation}
Following the same scheme as for Assertion 1, and under the
assumption that the volume of $M$ is infinite, we can prove that the
bottom of the spectrum of $\Delta_g$ in $M \setminus B(R_1)$, with
Dirichlet boundary conditions on $\partial B(R_1)$, satisfies the
inequality
$$\lambda_{\sigma}(\Delta_g, M\setminus B(R_1)) \le \dfrac{a}{4b-1}.$$

To conclude, we use the fact that
$$\lambda_{e}(\Delta_g) = \lim_{R \to \infty} \lambda_{\sigma}(\Delta_g, M\setminus B(R)).$$

\qed \bigskip



\begin{proposition}\label{P-vge1}
Let $(M,g)$ be a complete Riemannian surface with (non-negative)
Laplace operator $\Delta_g$ and Gaussian curvature $K_g$. Let $V(r)$
denote the volume of the geodesic ball of radius $r$ in $M$ (with
center some given point $x_0$). Let $a, b$ be positive real numbers,
with $b > 1/4$. Let $\alpha_0 = \sqrt{a/(4b-1)}$. If the operator $L
:= \Delta_g + a + b K_g$ has finite index, then
$$\forall \alpha > \alpha_0, ~~~ \int_0^{\infty} e^{-2\alpha r} V(r) \, dr < \infty
,$$ and hence, the lower volume growth of $M$ satisfies
$$\liminf_{r \to \infty} r^{-1} \ln (V(r))) \le 2 \alpha_0.$$
\end{proposition}

\pf It follows from our assumptions that the operator $L$ is
positive outside some compact set (see \cite{FC85}, Proposition 1).
In particular, it is positive on $M \setminus B(R_0)$ for some
radius $R_0$. Choose $R > R_0 + 1$ and define the function
\begin{equation}\label{E-vge-a1}
\xi (r) = \left\{%
\begin{array}{ll}
0 & \text{ ~ for ~ } r \le R_0,\\[4pt]
(1 - \frac{R_0+1}{R})^{\alpha R}(r-R_0) & \text{ ~ for ~ } R_0 \le r
\le R_0 + 1,\\[4pt]
(1 - \frac{r}{R})^{\alpha R} & \text{ ~ for ~ } R_0 + 1\le r \le R,
\end{array}%
\right.
\end{equation}
where the parameter $\alpha$ will be chosen later on. The positivity
of the operator $L$ on $M \setminus B(R_0)$ implies that
$$
0 \le \int_M \Big( (\xi'(r))^2 + a \xi^2(r) + b K_g \xi^2(r) \Big)
\, dv_g.
$$

We write the integral on the right-hand side as the sum of two
integrals, $\int_{C(R_0, R_0+1)}$ and $\int_{C(R_0+1, R)}$. The
first integral can be written as
$$
\int_{C(R_0, R_0+1)} = \Big( 1 - \frac{R_0+1}{R} \Big)^{\alpha R}
C(B(R_0)),
$$
where $C(B(R_0))$ is a constant which only depends on the geometry
of $M$ on the ball $B(R_0)$. Using Lemma \ref{L-cas} and the fact
that $\chi (B(r)) \le 1$ for all $r$, the second integral can be
estimated as follows
$$
\begin{array}{ll}
\displaystyle \int_{C(R_0+1, R)} \le & \displaystyle \int_{C(R_0+1,
R)} \Big( (\xi')^2 + a \xi^2 - b (\xi^2)''
\Big) \, dv_g \\[4pt]
& \\ & + 2\pi b - \xi^2(R_0+1)G(R_0+1) + 2\alpha L(R_0+1).
\end{array}
$$
Using (\ref{E-vge-a1}), the definition for the function $\xi$, the
integral in the first line of the above inequality can be written as
$$
- \Big( (4b-1) \alpha^2 - \frac{2b\alpha}{R} - a \Big)
\int_{R_0+1}^{R} \big( 1 - \frac{r}{R}\big)^{2\alpha R-2}L(r) \, dr.
$$

Taking $\alpha$ big enough so that the constant is positive, and
using the fact that $L(r) = V'(r)$, we obtain the inequality
$$
\frac{2\alpha R-2}{R} \Big( (4b-1) \alpha^2 - \frac{2b\alpha}{R} - a
\Big) \int_{R_0+1}^{R} \Big( 1 - \frac{r}{R}\Big)^{2\alpha R-3}V(r)
\, dr \le D(B(R_0),\alpha ),
$$
where $D(B(R_0),\alpha )$ is a constant which only depends on the
geometry of $M$ in the ball $B(R_0)$ and $\alpha$. Letting $R$ tend
to infinity, we finally obtain that

$$2\alpha \Big( (4b-1) \alpha^2 - a \Big) \int_{R_0+1}^{\infty}
e^{-2\alpha r}V(r) \, dr < \infty ,$$
provided that $\alpha > \alpha_0$, which proves the first assertion
in the theorem. The second assertion follows easily. \qed \bigskip

\textbf{Remark}. In the above theorem, we have assumed that $a>0$.
In the case $a=0$, one can show that the volume growth is at most
quadratic (see \cite{Cas06}, Proposition 2.2). \bigskip

\subsection{Applications to stable minimal surfaces in $\HH^3$ or  $\HH^2
\times \R$.}\label{ss-aps}

Let $M$ be a complete, orientable, minimal immersion into either the
$3$-dimensional hyperbolic space $\HH^3$ or into $\HH^2 \times \R$.
Let $J_M$ denote the Jacobi operator of the immersion. \bigskip

In the case of a minimal immersion $M \looparrowright \HH^3(-1)$,
the operator $J_M$ takes the form $J_M = \Delta_M + 2 - |A|^2$,
where $A$ is the second fundamental form. Using the Gauss equation,
we have that $K_M = - 1 - \frac{1}{2}|A|^2$, so that we can rewrite
the Jacobi operator of $M \looparrowright \HH^3(-1)$ as
\begin{equation}\label{E-jh3}
    J_M = \Delta_M + 4 + 2 K_M.
\end{equation}

In the case of a minimal immersion $M \looparrowright \HH^2(-1)
\times \R$, the Jacobi operator is given by $J_M = \Delta_M + 1 -
v^2 - |A|^2$, where $v$ is the vertical component of the unit normal
vector to the surface. Using the Gauss equation, we have that $K_M =
- v^2 - \frac{1}{2}|A|^2$, so that we can rewrite the Jacobi
operator of $M \looparrowright \HH^2(-1) \times \R$ as
\begin{equation}\label{E-jhr}
J_M = \Delta_M + 2 + 2K_M - (1-v^2) \le \tilde{J}_M := \Delta_M + 2
+ 2K_M.
\end{equation}

In this case, the positivity of the operator $J_M$ implies the
positivity of the operator $\tilde{J}_M$. \bigskip

Applying Proposition \ref{P-up} to the operator $J_M$ in the form
(\ref{E-jh3}) when $M$ is a minimal surface in $\HH^3$, \resp to the
operator $\tilde{J}_M$ in the form (\ref{E-jhr}) when $M$ is a
minimal surface in $\HH^2 \times \R$, we obtain the following
proposition.

\begin{prop}\label{P-vge-a2}
Let $(M,g) \looparrowright (\widehat{M}, \gh)$ be a complete,
orientable, minimal immersion. Assume that the immersion is stable.
\begin{enumerate}
    \item If $\widehat{M} = \HH^3$, then $\lambda_{\sigma}(\Delta_g) \le \frac{4}{7}.$
    \item If $\widehat{M} = \HH^2 \times \R$, then $\lambda_{\sigma}(\Delta_g)
    \le \frac{2}{7}.$
\end{enumerate}
If the immersion is only assumed to have finite index, then the same
inequalities hold with $\lambda_{\sigma}(\Delta_g)$ replaced by
$\lambda_{e}(\Delta_g)$, the infimum of the essential spectrum.
\end{prop}

\textbf{Remarks}. (i) The first assertion improves an earlier result
of A. Candel \cite{C} who proved that $\lambda_{\sigma}(M) \le
\frac{4}{3}$, provided that $M$ is a complete,
\emph{simply-connected}, stable minimal surface in $\HH^3$. (ii)
Note that in both cases, the bottom of the spectrum of a totally
geodesic $\HH^2$ is $1/4$.\bigskip


Applying Proposition \ref{P-vge1}, we have the following proposition.

\begin{prop}\label{P-vge-a3}
Let $(M,g) \looparrowright (\widehat{M}, \gh)$ be a complete,
orientable, minimal immersion. Let $\mu$ denote the lower volume
growth rate of $M$,
$$\mu = \liminf_{r \to \infty} r^{-1} \ln (V(r))),$$
where $V(r)$ is the volume of the geodesic ball $B(x_0,r)$ for some
given point $x_0$. Assume that the immersion has finite index.
\begin{enumerate}
    \item If $\widehat{M} = \HH^3$, then $\mu \le 2\sqrt{\frac{4}{7}}.$
    \item If $\widehat{M} = \HH^2 \times \R$, then $\mu \le 2\sqrt{\frac{2}{7}}.$
\end{enumerate}
\end{prop}\bigskip

\textbf{Remarks}. (i) Assertion 1 in Proposition \ref{P-vge-a3}
improves a previous result in \cite{C}, where Candel gives an upper
bound on $\mu$ under the assumption that $M$ is simply-connected.
(ii) Recall from \cite{K,B} that the volume growth is related to the
infimum of the essential spectrum by the formula

$$\lambda_{e}(\Delta_g) \le \Big( \frac{\liminf_{r \to \infty} r^{-1}
\ln (V(r)))}{2} \Big)^2.$$


\subsection{Futher applications}\label{ss-fa}

We note that the above argument also works for surfaces with
constant mean curvature  $|H| \leq 1$ in hyperbolic space. In that
case, $K_M=-(1-H^2)-\frac{1}{2}|A|^2$ and $J_M=\Delta_M + 4(1-H^2) +
2K_M$. So that, we obtain the following proposition. \smallskip

\begin{prop}\label{P-fa1}
Let $(M,g) \looparrowright \HH^3$ be a complete,
orientable, stable CMC  immersion, with $|H|\leq 1$. Then%
  $$
    \lambda_\sigma(\Delta_g)\leq \frac{4(1-H^2)}{7}.
  $$
\end{prop} \smallskip

The space $\HH^2\times\real$ is a simply-connected 3-dimensional
homogeneous manifold, whose isometry group has dimension $4$. Such
manifolds have been well studied  (see for instance \cite{D07} and
references therein) and can be parametrized by two real parameters,
say $\ka$ and $\tau$, with $\ka\neq 4\tau^2$. We denote them by
$\mathbb{E}^3(\kappa, \tau)$. When $\tau=0$, $\mathbb{E}^3(\kappa,
0)=$ is the pro\-duct space $\mathbb{E}^2(\ka)\times\real$, where
$\mathbb{E}^2(\ka)$ is the space form of constant curvature $\ka$.
In particular, $\HH^2\times\real = \mathbb{E}^3(-1,0)$.
\bigskip

If $(M,g) \looparrowright \mathbb{E}^3(\kappa, \tau)$ is an immersed
CMC $H$ surface, then its Jacobi operator is given by (see
\cite{D07}, Proposition 5.11)
$$
 J_M:=\Delta_g+2K-4H^2 -\ka - (\ka-4\tau^2)v^2.
$$

In the next proposition we give an upper bound for the bottom of the
spectrum in this general framework.

\begin{prop}\label{P-fa2}
Let $(M,g) \looparrowright \mathbb{E}^3(\kappa, \tau)$ be a
complete, orientable, stable CMC $H$ immersion, such that
$\ka<4\tau^2$. Assume furthermore that $2H^2\leq (2\tau^2 - \ka)$.
Then
$$
\la_\si(\De_g)\leq\frac{4 \tau^2 - 2\ka - 4H^2}{7}.
$$
\end{prop}
\bigskip

\pf Under the hypotheses we have the follows inequalities:
$$
0\leq \Delta_g+2K-4H^2 -\ka - (\ka-4\tau^2)v^2 \leq \Delta_g+2K-4H^2 - 2(\ka-2\tau^2),
$$
and we may apply Proposition \ref{P-up}  again. \qed

\section{Applications in higher dimensions}\label{s-apf}

In this Section, we give some further applications of the
inequalities we proved in Section \ref{s-hr}. In the following
proposition, we give a structure theorem for minimal hypersurfaces
in $\HH^m \times \R$.

\begin{proposition}\label{P-hr5}
Let $M^m \looparrowright \HH^m \times \R$, $m \ge 3$, be a complete,
orientable minimal hypersurface, with unit normal field $\nu$ and
second fundamental form $A$. Let $v$ denote the component of $\nu$
along $\partial_t$. For $0 \le \alpha \le 1$, there exists a
constant $c(m,\alpha)$ satisfying $c(m,\alpha) > 0$, whenever
\begin{enumerate}
    \item $m \ge 7$ and $\alpha \ge 0$,
    \item $m = 6$ and $\alpha \ge 0.083$,
    \item $m = 5$ and $\alpha \ge 0.578$.
\end{enumerate}
If $M$ satisfies $\|A\|_m \le c(m,\alpha)$ and $v^2 \ge \alpha^2$,
then $M$ carries no $L^2$-harmonic $1$-form and hence has at most
one end.
\end{proposition}\bigskip

\textbf{Proof.} We only sketch the proof. The proof uses several
ingredients. \bigskip

\textbf{1}. According to \cite{HS74}, the manifold $M^m$ satisfies
the Sobolev inequality
\begin{equation}\label{E-sob}
\|\varphi\|^2_{\frac{2m}{m-2}} \le S(2,m) \|d\varphi\|_2^2, ~~
\forall \varphi \in C_0^1(M).
\end{equation}

\textbf{2}. Let $u \in T_1M$ be a unit tangent vector to $M$. By
Gauss equation, we have the relation
$$\mathrm{Ric}(u,u) = \widehat{\mathrm{Ric}}(u,u) - \widehat{R}(u, \nu, u,
\nu) - |A(u)|^2,$$

where $\mathrm{Ric}$ denotes the Ricci curvature of $M$,
$\widehat{\mathrm{Ric}}$ the Ricci curvature and $\widehat{R}$ the
curvature tensor of $\Mh = \HH^m \times \R$, and where $A$ denotes
the Weingarten operator of the immersion. Using the curvature
computations in \cite{BS08} and the fact that $A$ has trace zero, we
obtain the inequality
\begin{equation}\label{E-apf1}
\mathrm{Ric}(u,u) \ge - (m-1) - \frac{m-1}{m}|A|^2.
\end{equation}

Let $\omega$ be an $L^2$ harmonic $1$-form on $M$. Using the
Weitzenb\"{o}ck formula for $1$-forms, the improved Kato inequality
\begin{equation}\label{E-kato}
\frac{1}{m-1} |d|\omega ||^2 \le |D\omega|^2 - |d|\omega ||^2,
\end{equation}

and inequality (\ref{E-apf1}), we find that $\omega$ satisfies the
following inequality in the weak sense
\begin{equation}\label{E-apf2}
\frac{1}{m-1} |d|\omega ||^2  + |\omega| \Delta |\omega| \le (m-1)
|\omega|^2 + \frac{m-1}{m}|A|^2 |\omega|^2.
\end{equation}

The following formal calculation can easily be made rigorous by
using cut-off functions. Integrate (\ref{E-apf2}) over $M$ using
integration by parts and the notation $f := |\omega|$,
$$\frac{m}{m-1} \int_M |df|^2 \le (m-1) \int_M f^2 + \frac{m-1}{m} \int_M |A|^2 f^2.$$

Plug the assumption $|v| \ge \alpha$ and the inequality
(\ref{E-hr7}) into the preceding inequality. Use H\"{o}lder's inequality
to estimate the integral $\int_M |A|^2 f^2$ and the Sobolev
inequality (\ref{E-sob}). It follows that

$$[\frac{m}{m-1} - \frac{4(m-1)}{(m-2+\alpha)^2}]
\|f\|_{\frac{2m}{m-2}}^2 \le S(2,m) \frac{m-1}{m} \|A\|_m^2
\|f\|_{\frac{2m}{m-2}}^2$$

and we can conclude the proof with the constant

$$C(m,\alpha) = \dfrac{m}{(m-1)(S(2,m))} \,
\dfrac{m(m-2+\alpha)^2-4(m-1)^2}{(m-1)(m-2+\alpha)^2}.$$

\qed \bigskip


\begin{proposition}\label{P-hr6}
Let $M^m \looparrowright \HH^m \times \R$ be a complete, orientable
minimal hypersurface, with second fundamental form $A$. Assume that
$||A||_m < \infty$. Then
\begin{enumerate}
    \item $M^m$ has finite index, if \, $m \ge 3$,
    \item $M^m$ has only finitely many ends,  if\, $m \ge 7$.
\end{enumerate}
\end{proposition}\bigskip

\textbf{Proof.} Assertion 1 was proved in \cite{BS08}. To prove
Assertion~2, we can mimic the proof of Corollary \ref{C-hr4} to show
that the operator $L:=\De+\frac{\sqrt{m-1}}{2}|A|^2-(m-1)$ has
finite index, when $m\geq 7$. We then apply Theorem 1 of \cite{dCWX}
to conclude the proof. \qed \smallskip


\begin{proposition}\label{P-hr8}
Let $M^m \looparrowright \HH^m \times \R$, $m \ge 3$, be a complete,
orientable minimal hypersurface, with unit normal field $\nu$ and
second fundamental form $A$. Let $v$ denote the component of $\nu$
along $\partial_t$. If,
\begin{enumerate}
    \item $\|A\|_{\infty} \le (\frac{m-1}{2})^2$, or
    \item $\|A\|_{\infty} \le (\frac{m-2+\alpha}{2})^2$ and $v^2 \ge
    \alpha$, or
    \item $|A|^2 + (m-1)v^2 \le \frac{m^2}{4}$ on $M$.
\end{enumerate}
then the immersion $M$ is stable.
\end{proposition}\bigskip

\textbf{Proof.} Recall that the Jacobi operator $J_M$ of the
immersion $M$ is given by the formula
$$J_M = \Delta_g + (m-1)(1-v^2) - |A|^2.$$

Assertion 1 follows from Proposition \ref{P-hr3}. Assertions 2 and 3
follow from Corollary \ref{C-hr2}. \qed \bigskip

\textbf{Remark 1}. The second condition is not so interesting
because it implies that $v$ does not vanish. If $M$ is connected, we
may assume that $v > 0$ and then $M$ is stable because $v$ is a
Jacobi field, $J_M (v) = 0$. \bigskip

\textbf{Remark 2}. We can write the operator $J_M$ as
$$J_M = \Delta_g - (\frac{m-2}{2})^2 + \big[(\frac{m}{2})^2 - |A|^2 \big].$$

In view of the results \`{a} la Lieb or Li-Yau, one can show that if the
integral
$$\int_M \big[(\frac{m}{2})^2 - |A|^2\big]_{-}^{m/2}$$

is small enough, then $M$ is stable. \bigskip









\bigskip

\begin{small}

\begin{tabular}{l}
Pierre B\'{e}rard\\
Universit\'{e} Grenoble 1\\
Institut Fourier (\textsc{ujf-cnrs})\\
B.P. 74\\
38402 Saint Martin d'H\`{e}res Cedex\\
France\\
\verb+Pierre.Berard@ujf-grenoble.fr+\\
\end{tabular}\bigskip

\begin{tabular}{l}
Philippe Castillon\\
Universit\'{e} Montpellier II\\
D\'{e}partement des sciences math\'{e}matiques CC 51\\
I3M (\textsc{umr 5149})\\
34095 Montpellier Cedex 5\\
France\\
\verb+cast@math.univ-montp2.fr+\\
\end{tabular}\bigskip

\begin{tabular}{l}
Marcos Cavalcante\\
Universidade Federal de Alagoas\\
Instituto de Matem\'{a}tica\\
57072-900 Macei\'o-AL\\
Brazil\\
\verb+marcos.petrucio@pq.cnpq.br+\\
(currently visiting Institut Fourier)
\end{tabular}\

\end{small}



\begin{thebibliography}{999}

\bibitem{BGS85}
\newblock{W. Balmann, M. Gromov and V. Schr\"{o}der,}
\newblock{\it Manifolds of nonpositive curvature,}
Progress in Math., {\bf 61}, Birkh\"{a}user 1985.


\bibitem{BdCS97}
\newblock{P. B\'{e}rard, M. do Carmo and W. Santos,}
\newblock{\it The index of constant mean curvature surfaces in hyperbolic 3-space,}
Math. Z., {\bf 224} (1997), 313--326.

\bibitem{BS08}
\newblock{P. B\'{e}rard and R. Sa Earp,}
\newblock{\it Minimal hypersurfaces in $\HH^n \times \R$, total curvature and index,}
arXiv:0808.3838v3.


\bibitem{B}
\newblock{ R. Brooks,}
\newblock{\it A relation between growth and the spectrum of the Laplacian,}
Math. Z., {\bf 178} (1981), no. 4, 501--508.

\bibitem{C}
\newblock{A. Candel,}
\newblock {\it Eigenvalue estimates for minimal surfaces in hyperbolic space,}
Trans. Amer. Math. Soc., {\bf 359} (2007), 3567--3575.


\bibitem{dCWX}
\newblock{M. do Carmo, Q. Wang and C. Xia,}
\newblock {\it Complete submanifolds with bounded mean curvature in a Hadamard
manifold,} J. of Geom. Phys., {\bf  60} (2010), 142--154.


\bibitem{Cas97}
\newblock{Ph. Castillon,}
\newblock {\it Sur l'op\'{e}rateur de stabilit\'{e} des sous-vari\'{e}t\'{e}s \`{a} courbure
moyenne constante dans l'espace hyperbolique,} Manuscripta Math.,
{\bf 94} (1997), 385--400.


\bibitem{Cas06}
\newblock{Ph. Castillon,}
\newblock {\it An inverse spectral problem on surfaces},
\newblock {Comment. Math. Helv.}, {\bf 81} (2006), 271--286.

\bibitem{CM02}
\newblock{T. Colding and W. Minicozzi},
\newblock {\it Estimates for parametric elliptic integrands,}
\newblock {Internat. Math. Res. Notices,} {\bf 6} (2002), 291--297.


\bibitem{CG92}
\newblock{J. Choe and R. Gulliver,}
\newblock {\it Isoperimetric inequalities on minimal submanifolds of space forms,}
Manuscripta Math., {\bf 77} (1992), 169--189.

\bibitem{D07}
\newblock{B. Daniel}
\newblock{\it Isometric immersions into 3-dimensional homogeneous
manifolds},
Comment. Math. Helv., {\bf 82} (2007), 87--131.

\bibitem{Da89}
\newblock{E. B. Davies,}
\newblock{\it Heat kernels and spectral theory,}
Cambridge University Press 1989.


\bibitem{FC85}
\newblock{D. Fischer-Colbrie,}
\newblock{\it On complete minimal surfaces with finite Morse index,}
Inventiones Math., {\bf 82} (1985), 121--132.


\bibitem{Gri99}
\newblock{A. Grigory'an,}
\newblock{\it Analytic  and  geometric  background  of
recurrence and  non-explosion of  the  brownian motion  on
Riemannian manifolds,} Bull. Amer. Math. Soc., {\bf 36} (1999),
135--249.

\bibitem{HS74}
\newblock{D. Hoffman and J. Spruck,}
\newblock{\it Sobolev and isoperimetric inequalities for Riemannian submanifods,}
Comm. Pure Applied Math., {\bf 27} (1974), 715--727.

\bibitem{K}
\newblock{H. Kumura,}
\newblock{\it Infimum of the exponential volume growth and the bottom
of the essential spectrum of the Laplacian,} arXiv:0707.0185v1
(2007).

\bibitem{LS97}
\newblock{D. Levin and M. Solomyak,}
\newblock{\it The Rozenblum-Lieb-Cwickel inequality for Markov generators,}
Journal d'analyse math\'{e}matique, {\bf 71} (1997), 173--193.

\bibitem{P81}
\newblock{A.V. Pogorelov,}
\newblock{\it On the stability of minimal surfaces,}
Soviet Math. Dokl., {\bf 24} (1981), 293--295.


\bibitem{ST88}
\newblock{K. Shiohama and M. Tanaka,}
\newblock{\it An isoperimetric problem for infinitely connected complete
open surfaces,} Geometry of manifolds (Matsumoto, 1988), Perspect.
Math. {\bf 8}, Academic Press, Boston, MA, (1989), 317--343.


\bibitem{ST93}
\newblock{K. Shiohama and M. Tanaka,}
\newblock{\it  The length function of geodesic parallel circle,}
Progress in differential geometry, (K. Shiohama, ed.), Adv. Stud.
Pure Math. {\bf 22}, Math. Soc. Japan, Tokyo (1993), 299--308.

\bibitem{Ty87}
\newblock{J. Tysk,}
\newblock{\it Eigenvalue estimates with applications to minimal surfaces,}
Pacific Journal of Math., {\bf 128} (1987), 361--366.




\end{thebibliography}
\end{document}